\documentclass[titlepage,twoside,12pt]{article}
\usepackage{amssymb}
\usepackage{amsfonts}
\begin{document}

{\bf \Large Solving Problems in Scalar Algebras \\ \\ of Reduced Powers} \\

{\bf Elem\'{e}r E Rosinger} \\
Department of Mathematics \\
and Applied Mathematics \\
University of Pretoria \\
Pretoria \\
0002 South Africa \\
eerosinger@hotmail.com \\

{\bf Abstract} \\

Following our previous work, we suggest here a large class of {\it algebras of scalars} in
which {\it simultaneous and correlated computations} can be performed owing to the existence
of surjective algebra homomorphisms. This may replace the currently used traditional
computations in which only real or complex scalars are used, or occasionally, nonstandard ones.
The usual real, complex, or nonstandard scalars are included in the mentioned large class of
algebras. \\
Such simultaneous and correlated computations offer a depth of insight which has so far been
missed when only using the few traditional kind of scalars. \\

{\bf 1. A Large Class of Scalar Algebras of Reduced Powers} \\

The following large class of {\it algebras of scalars} can be obtained easily as {\it reduced
powers}, Rosinger. This reduced power construction, in its more general forms, is one of the
fundamental tools in Model Theory, see Hodges. Historically, it has been used in the 19th
century in a particular case, by the classical Cauchy-Bolzano construction of the set
$\mathbb{R}$ of real numbers from the set $\mathbb{Q}$ of rational ones. \\

Let $\Lambda$ be any infinite set, then the power $\mathbb{R}^\Lambda$ is in a natural way an
associative and commutative {\it algebra}. Namely, the elements $x \in \mathbb{R}^\Lambda$ can
be seen as mappings $x : \Lambda \longrightarrow \mathbb{R}$, and as such, they can be added
to, and multiplied with one another point-wise. In the same way, the elements $x \in
\mathbb{R}^\Lambda$ can be multiplied with scalars from $\mathbb{R}$. \\

The well known remarkable fact connected with such a power algebra $\mathbb{R}^\Lambda$ is
that there is a one-to-one correspondence between the {\it proper ideals} in it, and on the
other hand, the {\it filters} on the infinite set $\Lambda$, Rosinger. Indeed, this one-to-one
correspondence operates as follows \\

(1.1)~~~ $ \begin{array}{l}
                {\cal I} ~~\longmapsto~~ {\cal F}_{\cal I} ~=~
                            \{~ Z(x) ~~|~~ x \in {\cal I} ~\} \\ \\
                {\cal F} ~~\longmapsto~~ {\cal I}_{\cal F} ~=~
                  \{~ x \in \Lambda \longrightarrow \mathbb{R} ~~|~~ Z(x) \in {\cal F} ~\}
           \end{array} $ \\

where ${\cal I}$ is an ideal in $\mathbb{R}^\Lambda$, ${\cal F}$ is a filter on $\Lambda$,
while for $ x \in \Lambda \longrightarrow \mathbb{R}$, we denoted $Z(x) = \{ \lambda \in
\Lambda ~|~ x(\lambda) = 0 \}$, that is, the zero set of $x$. \\

Important properties of the one-to-one correspondence in (1.1) are as follows. Given two
ideals ${\cal I}, {\cal J}$ in $\mathbb{R}^\Lambda$, and two filters ${\cal F}, {\cal G}$ on
$\Lambda$, we have \\

(1.2)~~~ $ \begin{array}{l}
             {\cal I} ~\subseteq~ {\cal J} ~~\Longrightarrow~~
                         {\cal F}_{\cal I} ~\subseteq~ {\cal F}_{\cal J} \\ \\
             {\cal F} ~\subseteq~ {\cal G} ~~\Longrightarrow~~
                         {\cal I}_{\cal F} ~\subseteq~ {\cal I}_{\cal G}
            \end{array} $ \\

Furthermore, the correspondences in (1.1) are idempotent when iterated, namely \\

(1.3)~~~ $ \begin{array}{l}
              {\cal I} ~~\longmapsto~~ {\cal F}_{\cal I}
                      ~~\longmapsto~~ {\cal I}_{{\cal F}_{\cal I}} ~=~ {\cal I} \\ \\
               {\cal F} ~~\longmapsto~~ {\cal I}_{\cal F}
                      ~~\longmapsto~~ {\cal F}_{{\cal I}_{\cal F}} ~=~ {\cal F}
            \end{array} $ \\

It follows that every {\it reduced power algebra} \\

(1.4)~~~ $ A ~=~ \mathbb{R}^\Lambda / {\cal I} $ \\

where ${\cal I}$ is an ideal in $\mathbb{R}^\Lambda$, is of the form \\

(1.5)~~~ $ A ~=~ A_{\cal F} ~\stackrel{def}=~ \mathbb{R}^\Lambda / {\cal I}_{\cal F} $ \\

for a suitable unique filter ${\cal F}$ on $\Lambda$. \\

Needless to say, the {\it advantage} of the representation of reduced power algebras given in
(1.5) is in the fact that filters  ${\cal F}$ on $\Lambda$ are {\it simpler} mathematical
structures, than ideals ${\cal I}$ in $\mathbb{R}^\Lambda$. \\

We shall call $\Lambda$ the {\it index set} of the reduced power algebra $A_{\cal F} ~=~
\mathbb{R}^\Lambda / {\cal I}_{\cal F}$, while ${\cal F}$ will be called the {\it generating
filter} which, we recall, is a filter on that index set. \\

Obviously, we can try to relate various reduced power algebra $A_{\cal F} ~=~
\mathbb{R}^\Lambda / {\cal I}_{\cal F}$ according to the two corresponding parameters which
define them, namely, their index sets and their generating filters. We start here by relating
them with respect to the latter. \\

Namely, a direct consequence of the second implication in (1.2) is the following one. Given
two filters ${\cal F} \subseteq {\cal G}$ on $\Lambda$, we have the {\it surjective algebra
homomorphism} \\

(1.6)~~~ $ A_{\cal F} \ni x + {\cal I}_{\cal F} ~~\longmapsto~~
                                      x + {\cal I}_{\cal G} \in A_{\cal G} $ \\

This obviously means that the algebra $A_{\cal G}$ is {\it smaller} than the algebra $A_{\cal
F}$, the precise meaning of it being that \\

(1.6$^*$)~~~ $ A_{\cal G} ~~\mbox{and}~~ A_{\cal F} /
               ( {\cal I}_{\cal G} / {\cal I}_{\cal F} ) ~~\mbox{are isomorphic algebras} $ \\

which follows from the so called {\it third isomorphism theorem for rings}, a classical result
of undergraduate Algebra. \\

Here we note that in the particular case when the filter ${\cal F}$ on $\Lambda$ is generated
by a nonvoid subset $I \subseteq \Lambda$, that is, when we have \\

(1.7)~~~ ${\cal F} = \{~ J \subseteq \Lambda ~|~ J \supseteq I ~\} $ \\

then it follows easily that \\

(1.8)~~~ $ A_{\cal F} ~=~ \mathbb{R}^I $ \\

which means that we do {\it not} in fact have a reduced power algebra, but only a power
algebra. \\
For instance, in case $I$ is finite and has $n \geq 1$ elements, then $A_{\cal F} =
\mathbb{R}^n$ is in fact the usual n-dimensional Euclidean space. \\

Consequently, in order to avoid such a degenerate case of reduced power algebras, we have to
avoid the filters of the form (1.7). This can be done easily, since such filters are obviously
characterized by the property \\

(1.9)~~~ $ \bigcap_{\,J\, \in\, {\cal F}}~ J ~=~ I \neq \phi $ \\

It follows that we shall only be interested in filters  ${\cal F}$ on $\Lambda$ which satisfy
the condition \\

(1.10)~~~ $ \bigcap_{\,J\, \in\, {\cal F}}~ J ~=~ \phi $ \\

or equivalently \\

(1.11)~~~ $ \begin{array}{l}
                \forall~~ \lambda \in \Lambda ~: \\ \\
                \exists~~ J_\lambda \in {\cal F} ~: \\ \\
                ~~~ \lambda \notin J_\lambda
            \end{array} $ \\

which is further equivalent with \\

(1.12)~~~ $ \begin{array}{l}
                \forall~~ I \subset \Lambda,~~ I ~~\mbox{finite} ~: \\ \\
                ~~~ \Lambda \setminus I \in {\cal F}
            \end{array} $ \\

We recall now that the {\it Frech\'{e}t filter} on $\Lambda$ is given by \\

(1.13)~~~ $ {\cal F}re ( \Lambda ) ~=~ \{~ \Lambda \setminus I ~~|~~
                             I \subset \Lambda,~~ I ~~\mbox{finite} ~\} $ \\

In this way, condition (1.10) - which we shall ask from now on about all filters ${\cal F}$ on
$\Lambda$ - can be written equivalently as \\

(1.14)~~~ $ {\cal F}re ( \Lambda ) ~\subseteq~ {\cal F} $ \\

This in particular means that \\

(1.14$^*$)~~~ $ \begin{array}{l}
                    \forall~~ I \in {\cal F} ~: \\ \\
                    ~~~ I ~~\mbox{is infinite}
                \end{array} $ \\

Indeed, if we have a finite $I \in {\cal F}$, then $\Lambda \setminus I \in {\cal F}re (
\Lambda )$, hence (1.14) gives $\Lambda \setminus I \in {\cal F}$. But $I \cap ( \Lambda
\setminus I ) = \phi$, and one of the axioms of filters is contradicted. \\

In view of (1.6), it follows that all reduced power algebras considered from now on will be
{\it homomorphic images} of the reduced power algebra $A_{{\cal F}re ( \Lambda )}$, through
the surjective algebra homomorphisms \\

(1.15)~~~ $ A_{{\cal F}re ( \Lambda )} \ni x + {\cal I}_{{\cal F}re ( \Lambda )}
                             ~~\longmapsto~~ x + {\cal I}_{\cal F} \in A_{\cal F} $ \\

or in view of (1.6$^*$), we have the {\it isomorphic algebras} \\

(1.15$^*$)~~~ $ A_{\cal F},~~~ A_{{\cal F}re ( \Lambda )} /
                        ( {\cal I}_{\cal F} /{\cal I}_{{\cal F}re ( \Lambda )} ) $ \\

Let us note that the {\it nonstandard reals}~ $^*\mathbb{R}$ are a particular case of the
above reduced power algebras (1.4). Indeed, $^*\mathbb{R}$ can be defined by using {\it free
ultrafilters} ${\cal F}$ on $\Lambda$, that is, ultrafilters which satisfy (1.10), or
equivalently (1.14). \\

We note that the field of real numbers $\mathbb{R}$ can be embedded naturally in each of the
reduced power algebras (1.4), by the {\it injective algebra homomorphism} \\

(1.16)~~~ $ \mathbb{R} \ni \xi ~~\longmapsto~~ x_\xi + {\cal I} \in A $ \\

where $x_\xi ( \lambda ) = \xi$, for $\lambda \in \Lambda$. Indeed, if $x_\xi \in {\cal I}$
and $\xi \neq 0$, then the ideal ${\cal I}$ must contain $x_1$, which means that it is {\it
not} a proper ideal, thus contradicting the assumption on it. \\

For simplicity of notation, from now on we shall write $x_\xi = \xi$, for $\xi \in \mathbb{R}$,
thus (1.16) will take the form \\

(1.17)~~~ $ \mathbb{R} \ni \xi ~~\longmapsto~~ \xi + {\cal I} \in A $ \\

which in view of the injectivity of this mapping, we may further simplify to \\

(1.18)~~~ $ \mathbb{R} \ni \xi ~~\longmapsto~~ \xi \in A $ \\

There is also the issue to relate reduced power algebras corresponding to different {\it index
sets}. Namely, let $\Lambda \subseteq \Gamma$ be two index sets which, as always in this paper,
are assumed to be both infinite. Then we have the obvious {\it surjective algebra
homomorphism} \\

(1.19)~~~ $ \mathbb{R}^\Gamma \ni x ~~\longmapsto~~ x|_{\,\Lambda} \in \mathbb{R}^\Lambda $ \\

since the elements $x \in \mathbb{R}^\Gamma$ can be seen as mappings $x : \Gamma
\longrightarrow \mathbb{R}$. Consequently, given any ideal ${\cal I}$ in $\mathbb{R}^\Gamma$,
we can associate with it the ideal in $\mathbb{R}^\Lambda$, given by \\

(1.20)~~~ $ {\cal I}|_{\,\Lambda} ~=~ \{~ x|_{\,\Lambda} ~~|~~ x \in {\cal I} ~\} $ \\

As it happens, however, such an ideal ${\cal I}|_{\,\Lambda}$ need not always be a proper
ideal in $\mathbb{R}^\Lambda$, even if ${\cal I}$ is a proper ideal in $\mathbb{R}^\Gamma$.
For instance, if we take $\gamma \in \Gamma \setminus \Lambda$, and consider the proper ideal
in $\mathbb{R}^\Gamma$ given by ${\cal I} = \{ x \in \mathbb{R}^\Gamma ~|~ x(\gamma) = 0 ~\}$,
then we obtain ${\cal I}|_{\,\Lambda} = \mathbb{R}^\Lambda$, which is not a proper ideal in
$\mathbb{R}^\Lambda$. \\

We can avoid that difficulty by noting the following. Given a filter ${\cal F}$ on $\Gamma$
which satisfies (1.14), that is, ${\cal F}re ( \Gamma ) \subseteq~ {\cal F}$, then \\

(1.21)~~~ $ {\cal F}|_{\,\Lambda} ~=~ \{~ I \cap \Lambda ~~|~~ I \in {\cal F} ~\} $ \\

satisfies the corresponding version of (1.14), namely ${\cal F}re ( \Lambda ) ~\subseteq~
{\cal F}|_{\,\Lambda}$. Indeed, let us take $J \subseteq \Lambda$ such that $\Lambda \setminus
J$ is finite. Then clearly $\Gamma \setminus ( J \cup ( \Gamma \setminus \Lambda ) )$ is
finite, hence $J \cup ( \Gamma \setminus \Lambda ) \in {\cal F}$. However, $J = ( J \cup (
\Gamma \setminus \Lambda ) ) \cap \Lambda )$, thus $J \in {\cal F}|_{\,\Lambda}$. \\
Now in order for ${\cal F}|_{\,\Lambda}$ to be a filter on $\Lambda$, it suffices to show that
$\phi \notin {\cal F}|_{\,\Lambda}$. Assume on the contrary that for some $I \in {\cal F}$ we
have $I \cap \Lambda = \phi$, then $I \subseteq \Gamma \setminus \Lambda$, thus $\Lambda
\notin {\cal F}$. \\

It follows that \\

(1.22)~~~ $ {\cal F}|_{\,\Lambda} ~~\mbox{is a filter on}~ \Lambda
                         ~\mbox{which satisfies}~ (1.14)
                                    ~~\Longleftrightarrow~~ \Lambda \in {\cal F} $ \\

In view of (1.19) - (1.22), for every filter ${\cal F}$ on $\Gamma$, such that \\

(1.23)~~~ $ \Lambda \in {\cal F} $ \\

we obtain the {\it surjective algebra homomorphism} \\

(1.24)~~~ $ A_{\cal F} ~=~ \mathbb{R}^\Gamma / {\cal I}_{\cal F} \ni
               x + {\cal I}_{\cal F} ~~\longmapsto~~
                 x|_{\,\Lambda} + {\cal I}_{{\cal F}|_{\,\Lambda}} \in
                    A_{{\cal F}|_{\,\Lambda}} ~=~
                           \mathbb{R}^\Lambda / {\cal I}_{{\cal F}|_{\,\Lambda}} $ \\

and in particular, we have the following relation between the respective proper ideals \\

(1.25)~~~ $ ( {\cal I}_{\cal F} )|_{\,\Lambda} ~=~ {\cal I}_{{\cal F}|_{\,\Lambda}} $ \\

{\bf 2. Zero Divisors and the Archimedean Property} \\

It is an elementary fact of Algebra that a quotient algebra (1.4) has {\it zero divisors},
unless the ideal ${\cal I}$ is {\it prime}. A particular case of that is when a quotient
algebra (1.4) is a {\it field}, which is characterized by the ideal ${\cal I}$ being {\it
maximal}. And in view of (1.5), (1.2), this means that the filter ${\cal F}$ generating such
an ideal must be an {\it ultrafilter}. \\

On the other hand, {\it none} of the reduced power algebras (1.5) which correspond to filters
satisfying (1.14) are {\it Archimedean}. And that includes the nonstandard reals
$^*\mathbb{R}$ as well. \\

Let us elaborate in some detail in this regard. First we note that on reduced power algebras
(1.5), one can naturally define a {\it partial order} as follows. Given two elements $x +
{\cal I}_{\cal F},~ y + {\cal I}_{\cal F} \in A_{\cal F} = \mathbb{R}^\Lambda /
{\cal I}_{\cal F}$, we define \\

(2.1)~~~ $ x + {\cal I}_{\cal F} ~\leq~ y + {\cal I}_{\cal F} ~~\Longleftrightarrow~~
                    \{~ \lambda \in \Lambda ~~|~~
                           x ( \lambda ) ~\leq~ y ( \lambda ) ~\} \in {\cal F} $ \\

Now, with this partial order, the algebra $A_{\cal F}$ is called {\it Archimedean}, if and
only if \\

(2.2)~~~ $ \begin{array}{l}
               \exists~~ u + {\cal I}_{\cal F} \in A_{\cal F},~
                                            u + {\cal I}_{\cal F} \geq 0 ~: \\ \\
               \forall~~ x + {\cal I}_{\cal F} \in A_{\cal F},~
                                            x + {\cal I}_{\cal F} \geq 0 ~: \\ \\
               \exists~~ n \in \mathbb{N} ~: \\ \\
               ~~~  x + {\cal I}_{\cal F} ~\leq~ n (  u + {\cal I}_{\cal F} )
            \end{array} $ \\

However, in view of (1.14), we can take an infinite $I \in {\cal F}$. Thus we can define a
mapping $v : \Lambda \longrightarrow \mathbb{R}$ which is unbounded for above on $I$. And in
this case taking $x + {\cal I}_{\cal F} = ( u + v ) + {\cal I}_{\cal F}$, it follows easily
that condition (2.2) is not satisfied. \\

We note that the reduced power algebras (1.5) are Archimedean only in the degenerate case
(1.7), (1.8), when in addition the respective sets $I$ are {\it finite}, thus as noted, the
respective algebras reduce to finite dimensional Euclidean spaces. \\

{\bf 3. Simultaneous Correlated Computations in Scalar Algebras \\
        \hspace*{0.4cm} of Reduced Powers} \\

Instead of the traditional way which confines all computations with scalars to the real
numbers in $\mathbb{R}$ or to the complex numbers in $\mathbb{C}$, one can do a {\it
simultaneous and correlated} scalar computation in {\it all} of the reduced power algebras, by
using the surjective algebra homomorphisms (1.15) and (1.24). \\
Such a computation offers a depth of insight which has so far been missed when computing only
with scalars in $\mathbb{R}$, $\mathbb{C}$, or even in the nonstandard $^*\mathbb{R}$, all of
which are included as particular cases in the large, {\it two parameter} class of reduced
power algebras $A_{\cal F} = \mathbb{R}^\Lambda / {\cal I}_{\cal F}$, where $\Lambda$ can be
any {\it infinite set}, while ${\cal F}$ can be any {\it filter} on $\Lambda$ which satisfies
(1.14). \\

Let us consider in some more detail this two parameter family of reduced power algebras. \\

Let us start by fixing any given infinite index set $\Lambda$. Then, when ${\cal F}$ ranges
over all the filters on $\Lambda$ which satisfy (1.14), the {\it largest} corresponding
reduced power algebra of type (1.4) is given by, see (1.15) \\

(3.1)~~~ $ A_{{\cal F}re ( \Lambda )} ~=~
                 \mathbb{R}^\Lambda / {\cal I}_{{\cal F}re ( \Lambda )} $ \\

and all other reduced power algebras $A_{\cal F} = \mathbb{R}^\Lambda / {\cal F}$, with ${\cal
F}$ satisfying (1.14), are images of it through the surjective algebra homomorphism (1.15). \\

Now, let us allow the index set $\Lambda$ to range over all infinite sets. Then the algebras
(3.1) become ordered by the surjective algebra homomorphisms (1.24), namely \\

(3.2)~~~ $ A_{{\cal F}re ( \Gamma )} ~~\longrightarrow~~ A_{{\cal F}re ( \Lambda )} $ \\

whenever two infinite index sets $\Lambda$ and $\Gamma$ are in the relation \\

(3.3)~~~ $ \Lambda \subseteq \Gamma,~~~ \Gamma \setminus \Lambda ~~\mbox{is finite} $ \\

in which case we clearly have \\

(3.4)~~~ $ {{\cal F}re ( \Lambda )} ~=~ {{\cal F}re ( \Gamma )}|_{\,\Lambda} $ \\

Indeed, for $\Lambda \subseteq \Gamma$, we obviously have, see (1.23) \\

(3.5)~~~ $ \Lambda \in {\cal F}re ( \Gamma ) ~~\Longleftrightarrow~~
                   \Gamma \setminus \Lambda ~~\mbox{is finite} $ \\

thus according to (1.23), (1.24), the surjective algebra homomorphism (3.2) holds. \\

In this way we obtain the following {\it two directional} table of {\it surjective algebra
homomorphisms} \\

(3.6)~~~ $ \begin{array}{l}
               .~.~.~.~.~.~.~.~.~.~.~.~. \\
               \,\downarrow~~~~~~~~~~~~~~~~~~~~~~~~~~~~\downarrow~~~~~~~~~~~~~~~~~~~~~~
                                                                    \downarrow \\
               A_{{\cal F}re ( \Gamma )} ~\longrightarrow ~.~.~.~ \longrightarrow~ A_{\cal F}
                            ~\longrightarrow
                ~.~.~.~ \longrightarrow A_{\cal G} ~\longrightarrow ~.~.~.~ \\
               \,\downarrow~~~~~~~~~~~~~~~~~~~~~~~~~~~~\downarrow~~~~~~~~~~~~~~~~~~~~~~
                                                                    \downarrow \\
               .~.~.~.~.~.~.~.~.~.~.~.~. \\
               \,\downarrow~~~~~~~~~~~~~~~~~~~~~~~~~~~~\downarrow~~~~~~~~~~~~~~~~~~~~~~
                                                                    \downarrow \\
               A_{{\cal F}re ( \Lambda )} ~\longrightarrow~.~.~.~ \longrightarrow~ A_{\cal H}
                             ~\longrightarrow
                ~.~.~.~ \longrightarrow A_{\cal K} ~\longrightarrow ~.~.~.~ \\
               \,\downarrow~~~~~~~~~~~~~~~~~~~~~~~~~~~~\downarrow~~~~~~~~~~~~~~~~~~~~~~
                                                                    \downarrow \\
               .~.~.~.~.~.~.~.~.~.~.~.~.
            \end{array} $ \\

whenever the following conditions are satisfied \\

(3.6$^*$)~~~ $ \begin{array}{l}
              \Gamma ~\supseteq~ \Lambda,~~~ \Gamma \setminus \Lambda ~~\mbox{is finite} \\ \\
              {\cal F},~ {\cal G} ~~\mbox{are filters on}~ \Gamma,~~
                          {\cal F}re ( \Gamma ) \subseteq {\cal F} \subseteq {\cal G} \\ \\
              {\cal H},~ {\cal K} ~~\mbox{are filters on}~ \Lambda,~~
                          {\cal F}re ( \Lambda ) \subseteq {\cal H} \subseteq {\cal K} \\ \\
              {\cal F}|_{\,\Lambda} \subseteq {\cal H},~~ {\cal G}|_{\,\Lambda} \subseteq
                                                                                   {\cal K}
                \end{array} $ \\

Furthermore, as follows easily from (1.15), (1.24), under the above conditions (3.6$^*$), the
following diagrams of surjective algebra homomorphisms are {\it commutative} \\

(3.7)~~~ $ \begin{array}{l}
           A_{{\cal F}re ( \Gamma )} ~\longrightarrow~ A_{\cal F} ~\longrightarrow~
                                                  A_{\cal G} \\
           \,\downarrow~~~~~~~~~~~~~~~\downarrow~~~~~~~~~~\downarrow \\

               A_{{\cal F}re ( \Lambda )}~\longrightarrow~ A_{\cal H} ~\longrightarrow~
                                                  A_{\cal K}
            \end{array} $ \\

Finally, if we disregard the leftmost column in (3.6), then under {\it weaker} conditions than
in (3.6$^*$), we still obtain the following {\it commutative diagrams of surjective algebra
homomorphisms} \\

(3.8)~~~ $ \begin{array}{l}
               .~.~.~.~.~.~.~.~.~.~.~.~. \\
               ~~~~~~~~~~~~~~\downarrow~~~~~~~~~~~~~~~~~~~~~~\downarrow \\
                ~.~.~.~ \longrightarrow~ A_{\cal F} ~\longrightarrow
                ~.~.~.~ \longrightarrow A_{\cal G} ~\longrightarrow ~.~.~.~ \\
               ~~~~~~~~~~~~~~\downarrow~~~~~~~~~~~~~~~~~~~~~~\downarrow \\
               .~.~.~.~.~.~.~.~.~.~.~.~. \\
               ~~~~~~~~~~~~~~\downarrow~~~~~~~~~~~~~~~~~~~~~~\downarrow \\
               ~.~.~.~ \longrightarrow~ A_{\cal H}~\longrightarrow
                ~.~.~.~ \longrightarrow A_{\cal K} ~\longrightarrow ~.~.~.~ \\
               ~~~~~~~~~~~~~~\downarrow~~~~~~~~~~~~~~~~~~~~~~\downarrow \\
               .~.~.~.~.~.~.~.~.~.~.~.~.
            \end{array} $ \\

which hold whenever \\

(3.8$^*$)~~~ $ \begin{array}{l}
              \Gamma ~\supseteq~ \Lambda \\ \\
              {\cal F},~ {\cal G} ~~\mbox{are filters on}~ \Gamma,~~
                         {\cal F}re ( \Gamma ) \subseteq {\cal F} \subseteq {\cal G},~~
                                       \Lambda \in {\cal F} \\ \\
              {\cal H},~ {\cal K} ~~\mbox{are filters on}~ \Lambda,~~
                          {\cal F}re ( \Lambda ) \subseteq {\cal H} \subseteq {\cal K} \\ \\
              {\cal F}|_{\,\Lambda} \subseteq {\cal H},~~ {\cal G}|_{\,\Lambda} \subseteq
                                                                                   {\cal K}
                \end{array} $ \\

In this case, instead of (3.7), we obtain the following {\it commutative} diagrams of
surjective algebra homomorphisms \\

(3.9)~~~ $ \begin{array}{l}
                     A_{\cal F} ~\longrightarrow~ A_{\cal G} \\
                     \,\downarrow~~~~~~~~~~\downarrow \\
                     A_{\cal H} ~\longrightarrow~ A_{\cal K}
           \end{array} $ \\

In subsequent papers, we shall apply the above method of {\it simultaneous and correlated}
scalar computation to several important problems in Theoretical Physics. \\ \\

\end{document}